\documentclass[12pt]{amsart}

\usepackage[all]{xy}

\CompileMatrices

\theoremstyle{plain}
\newtheorem{lemma}{Lemma}[section]
\newtheorem{teo}[lemma]{Theorem}
\newtheorem{propo}[lemma]{Proposition}

\theoremstyle{definition}
\newtheorem{defi}[lemma]{Definition}

\theoremstyle{remark} 
\newtheorem{remark}[lemma]{Remark}
\newtheorem{con}[lemma]{Conjecture}

\newcommand{\pic}{{\rm Pic}}
\newcommand{\bs}{{\rm Bs}}
\newcommand{\p}{\mathbb{P}}
\newcommand{\co}{\mathbb{C}}
\newcommand{\z}{\mathbb{Z}}
\newcommand{\f}{\mathbb{F}}
\newcommand{\n}{\mathbb{N}}
\newcommand{\oc}{{\mathcal O}}
\newcommand{\ls}{{\mathcal L}}
\newcommand{\ms}{{\mathcal M}}
\newcommand{\rs}{{\mathcal R}}
\newcommand{\lss}{{\mathcal S}}
\newcommand{\lfe}{\ls_{\tilde{\f}}}
\newcommand{\lfn}{\ls_{\f}}
\newcommand{\rfn}{\rs_{\tilde{\f}}}
\newcommand{\re}{R_{\f}}
\newcommand{\rn}{R_{\tilde{\f}}}
\newcommand{\rfe}{\rs_{\f}}
\newcommand{\drfe}{r_{\tilde{\f}}}
\newcommand{\drfn}{r_{\f}}
\newcommand{\klfe}{\hat{\ls}_{\tilde{\f}}}
\newcommand{\klfn}{\hat{\ls}_{\f}}
\newcommand{\vfe}{v_{\tilde{\f}}}
\newcommand{\vfn}{v_{\f}}
\newcommand{\kvfe}{\hat{v}_{\tilde{\f}}}
\newcommand{\kvfn}{\hat{v}_{\f}}
\newcommand{\dlfe}{l_{\tilde{\f}}}
\newcommand{\dlfn}{l_{\f}}
\newcommand{\dklfe}{\hat{l}_{\tilde{\f}}}
\newcommand{\dklfn}{\hat{l}_{\f}}
\newcommand{\rt}{\rightarrow}

\begin{document}
\setlength{\baselineskip}{15pt}

\title{On linear systems of curves on rational scrolls}
\author{Antonio Laface}
\email{laface@mat.unimi.it}
\address{Dipartimento di Matematica, Universit\`a degli Studi di Milano 
Via Saldini, 50 20133 MILANO}
\begin{abstract}
In this paper we prove a conjecture on the dimension of linear systems, 
with base points of multiplicity $2$ and $3$, on an Hirzebruck surface.
\end{abstract}
\maketitle

\section{Introduction}

Consider $r$ points in general position on an algebraic surface $S$, 
to each one of these points $p_i$ we associate a natural number $m_i$ 
called the {\em multiplicity} of the point. 
Let $r_j$ be the number of $p_i$ with multiplicity $m_i$ and let $\ls$ be the 
complete linear system associated to the line bundle $L\in\pic(S)$
By $\ls({m_1}^{r_1},\cdots ,{m_k}^{r_k})$ we mean the linear system 
of curves in $\ls$ with $r_j$ general base points of
multiplicity at least $m_j$ for $j=1\cdots k$. 
Now, define the {\em effective dimension} of the system to be:
\[
l({m_1}^{r_1},\cdots ,{m_k}^{r_k})=\dim 
\ls({m_1}^{r_1},\cdots ,{m_k}^{r_k}),
\]
and the {\em virtual dimension} to be:
\[
v({m_1}^{r_1},\cdots ,{m_k}^{r_k})=
\dim\ls-\sum r_i\frac{m_i(m_i+1)}{2}.
\]
Observe that it may happen that $v<-1$, in
this case $v$ does not represent the dimension of a linear system and
instead of it we may define the {\em expected dimension}
$e=\max\{v,-1\}$.
From the previous definitions it follows immediately that for a given
system we have:
\begin{equation}
\label{remark1}
v\leq e\leq l.
\end{equation}

Let $Z$ be the $0$-dimensional scheme defined by the multiple points $p_i$
and consider the exact sequence of sheaves:
\[
0\rt {\mathcal I}(Z)\rt\oc_{S}\rt\oc_Z\rt 0
\]
where ${\mathcal I}(Z)$ is the ideal sheaf of $Z$. 
Tensoring with $L$ and taking cohomology we obtain:
\[
0\rt H^0(L\otimes {\mathcal I}(Z))\rt H^0(L)\rt H^0(L_Z)\rt 
\]
\[
\rt H^1(L\otimes {\mathcal I}(Z))\rt H^1(L)\rt 0.
\]
In this way we see that 
\begin{equation}
\label{equx}
h^1(L\otimes {\mathcal I}(Z))-h^1(L)=
l({m_1}^{r_1},\cdots ,{m_k}^{r_k}) - 
v({m_1}^{r_1},\cdots ,{m_k}^{r_k})
\end{equation}
Observe that the second inequality of (\ref{remark1}) may be strict since
the conditions imposed by the points may fail to be independent. 
In this case we say that the system is {\em special}.
The aim of this paper is to give a characterization of special systems
on rational scrolls and a complete classification will be given in the
case of homogeneous systems of multiplicity $\leq 3$. \\
The paper is organized as follows:\\
In section 2 we give some preliminaries on rational scrolls and linear
systems on them. Then we give an example of special systems $\ls$ whose
base locus contains a particular kind of multiple curves which  
we will call $(-1)$-curves. This will allow us to make a general conjecture 
on the structure of special linear systems.
In section 3 we restrict our attention to linear systems with point
of equal multiplicity, which we call homogeneous. Then we give a 
classification of homogeneous $(-1)$-curves of multiplicity $1$. 
In section 4 we study a class of birational 
transformations from $\f_n\rt\f_{n-1}$ which are a natural generalization of
quadratic transformations on $\p^2$. 
In section 5 we give the complete list of homogeneous $(-1)$-special
systems with multiplicity $m\leq 3$. 
In section 6 we describe the fundamental technique, which consists 
in a degeneration of the $\f_n$ to a suitable reducible surface.
In section 7 we prove the conjecture for all homogeneous systems 
systems with base points of multiplicity $\leq 3$.

\section{Preliminaries}

Consider an $\f_n$ surface, i.e. $\f_n=\p (\oc_{\p^1}\oplus\oc_{\p^1}(n))$
with $n\in\n$ and let $F,H$ be the two
generators of $\pic (\f_n)$ such that $F^2=0,\ H^2=n,\ F\cdot H=1$. 
With $\Gamma_n$ we indicate the {\em $(-n)$-curve} of $\f_n$ (i.e. the rational
curve $\Gamma_n\in\mid H-nF\mid$ of self intersection $-n$).
Let $L=aF+bH$ be a divisor, with $\bs\mid L\mid$ we mean the base locus of
the linear system $\mid L\mid$.
\begin{propo}
On an $\f_n$ surface, the ample divisors coincide with the very ample ones 
and are of the form $aF+bH$ with $a,b>0$. 
\end{propo}

\begin{proof}
\cite[Cap. V, 2.18]{ha}: 
\end{proof}
\begin{propo}
\label{effective}
Let $L=aF+bH$ be an effective divisor on an $\f_n$ surface, then it is linearly
equivalent to a divisor of the form: $q\Gamma_n +a_1F+b_1H$, with
$\bs\mid L\mid=q\Gamma_n$ and $a_1,b_1\geq 0$.
\end{propo}

\begin{proof}
Consider the product $L\cdot F=b$, this
must be $\geq 0$ since otherwise $F$ would be a fixed component of $\mid
L\mid$, which is impossible. Now consider the product 
$L\cdot\Gamma_n=a$, if $a<0$ then $\Gamma_n$ is contained in the fixed part of  
$\mid L\mid$. Now find $q,r\in\n$ such that $-a=qn-r$ with $r<n$. Observe
that $(L-(q-1)E)\cdot\Gamma_n<0$, hence $q\Gamma_n\in\bs\mid L\mid$.
Now consider the system: $\mid L-q\Gamma_n\mid$.
This is equal to: $\mid aF+bH-q\Gamma_n\mid = \mid aF+bH+qnF-qH\mid =\mid
rF+(b-q)H\mid = \mid a_1F+b_1H\mid$. 
From $L\cdot H\geq 0$ we obtain $bn+a\geq 0$ which gives $bn\geq qn-r$ and 
$r\geq (q-b)n$, so $b\geq q$ otherwise $r$ would be a positive multiple of
$n$. Hence both $a_1$ and $b_1$ are $\geq 0$. This implies that $\bs\mid
a_1F+b_1H\mid=\emptyset$. 
\end{proof}
In order to compute $h^0(aF+bH)$ (with $a,b\geq 0$) we need the following
proposition:
\begin{propo}
\label{vanish}
Let $aF+bH$ be a divisor on an $\f_n$ surface such that $a\geq 0$ and 
$b\geq -1$, then: 
\[
h^1(aF+bH)=0 \hspace{1cm} h^0(aF+bH)=(b+1)(2a+2+nb)/2.
\]
\end{propo} 

\begin{proof}
Consider the exact sequence:
\[
0\rt\oc_{\f_n}(-F) \rt\oc_{\f_n}\rt\oc_F\rt 0,
\]
tensoring with $\oc_{\f_n}(aF+bH)$ we obtain surjective maps for each
$a\in\z$:
\[
\cdots\rt H^1((a-1)F+bH)\rt H^1(aF+bH)\rt 0.
\]
Again from the sequence:
\[
0\rt\oc_{\f_n}(-H) \rt\oc_{\f_n}\rt\oc_H\rt 0,
\]
tensoring with $\oc_{\f_n}(bH)$, with $b\geq 0$, we obtain:
\[
\cdots\rt H^1((b-1)H)\rt H^1(bH)\rt 0.
\]
Now $h^1(-H)=0$ since $H^1(\oc_{\f_n})=0$ and the map $H^0(\oc_{\f_n})\rt
H^0(\oc_H)$ is an isomorphism, hence $h^1(aF+bH)=0$. \\
By Riemann-Roch and the preceding discussion, it follows that 
$h^0(aF+bH)=((aF+bH)^2-(aF+bH)\cdot K)/2+1=(b+1)(2a+2+nb)/2$. 
\end{proof}

Let $\ls_n(a,b,{m_1}^{r_1},\cdots ,{m_k}^{r_k})$ denote the linear system
$\mid aF+bH\mid$ on the $\f_n$ surface with $r_i$ base points 
of multiplicity $m_i$. 
By a generality assumption, we may suppose that no one of these points 
lies on the $(-n)$-curve $\Gamma_n$,
hence, by \ref{effective} we are lead to consider only systems 
with $a,b\geq 0$. \\

In what follows we will assume that our $r$ points on $\f_n$ will be in 
{\em general position}; this means that we choose the points in such a way that
each linear
systems with these simple base points in non-special. The scheme that
parameterize such points is a Zariski open set of the scheme that parameterize 
$r$ points on an $\f_n$ surface (this is why we use the word "generic"). \\

\begin{propo}
\label{mult.one}
A linear system $\ls$ on an $\f_n$ surface, with base points of multiplicity one
and in general position, is non-special.
\end{propo}

Under the assumption of general position we consider a particular class of 
special system, obtained in the following way: \\

Consider the blow up $S$ of $\f_n$ at $p_1,\ldots ,p_r$ and
denote by $\ls$ the strict transform of the system 
$\ls_n(a,b,{m_1}^{r_1},\cdots ,{m_k}^{r_k})$. 
Define the virtual and the expected dimension of $\ls$ as those of 
$\ls_n(a,b,{m_1}^{r_1},\cdots ,{m_k}^{r_k})$. 
By Riemann-Roch, the virtual dimension of $\ls$ may be given 
in the following way:
\[
v(\ls)=\frac{\ls^2-\ls\cdot K_S}{2},
\]
where $K_S$ denote the canonical bundle of $S$.
We recall that $E$ is a $(-1)$-curve on $S$ if $E$ is irreducible and 
$E^2=E\cdot K_S=-1$. 
Suppose that $\ls\cdot E=-t<0$, then $tE$ is contained in the fixed part of 
$\ls$. Let $\ms=\ls-tE$ be the residual system, then: 
\[
v(\ls)=\frac{\ls^2-\ls\cdot K_S}{2}=
\frac{(\ms+tE)^2-(\ms+tE)\cdot K_S}{2}=v(\ms)+\frac{t-t^2}{2}
\]
So if $t\geq 2$ and $\ls$ is not empty, then 
$l(\ls)=l(\ms)\geq v(\ms)>v(\ls)$, which means that the system $\ls$ is 
special. \\
If $t=1$, then $v(\ls)=v(\ms)$, but it may happen that $\ms\cdot\Gamma_n <0$.
In this case the $(-n)$-curve $\Gamma_n$ is contained in the fixed part of $\ls$
and the system may again be special. For example consider the system
$\ls_6(0,4,3^{11})$, this has expected dimension equal to $-1$.
Consider the $(-1)$-curve $E=\ls_6(2,1,1^{11})$, then $\ls\cdot E=-1$. 
The system $\ls-E=\ls_6(-2,3,2^{11})$ has again a fixed part, since it
has negative product with the $(-n)$-curve. The residual system
is $\ls-E-\Gamma_n= \ls_6(4,2,2^{11})$, and this is exactly 
$2E$, hence $\ls=\Gamma_n +3E$ and the initial system is not empty. \\

Given a linear system $\ls$, consider the following procedure:
\begin{itemize}
\item[1)] if it does exist a $(-1)$-curve $E$ such that
               $t:=\ls\cdot E<0$, then substitute $\ls$ with $\ls-tE$ and
               goto step 1), else goto step 2).
\item[2)] if $\ls\cdot\Gamma_n<0$ then substitute $\ls$ with
               $\ls-\Gamma_n$ and goto step 1), else finish.
\end{itemize}
After this procedure, that must obviously ends in a finite number of
steps, we have a new linear system $\ms$:

\begin{defi}
A linear system $\ls$ is \emph{$(-1)$-special} if $v(\ls)>v(\ms)$ 
\end{defi}

For linear systems on $\p^2$ a conjecture due to Hirschowitz says that
the only special systems are the $(-1)$-special ones. Reformulating this 
conjecture on $\f_n$ surfaces we have:

\begin{con}
A linear system $\ls_n(a,b,{m_1}^{r_1},\cdots ,{m_k}^{r_k})$ is special if and
only if it is $(-1)$-special.
\end{con}
If $m_1=\cdots m_k=m$ the system is called homogeneous of multiplicity $m$.
The aim of this paper is to prove the conjecture for homogeneous systems with
$m\leq 3$. 
This conjecture may be easily proved in the following case:
\begin{propo}
\label{propo:<=n+1}
The conjecture is true on an $\f_n$ surface if the number $r$ of points 
is less then or equal to $n+3$.
\end{propo}
\begin{proof}
It is enough to prove the conjecture when $r=n+3$, so from now 
on all the sums are intended to be over the $n+3$ points.

Let $\pi: S\to\f_n$ be the blow-up map and observe that $-K_S=N+\Delta$, where
$N:=\pi^*(\frac{2}{n}\Gamma_n+H+2F)-\sum_iE_i$ and $\Delta:=\pi^*(\frac{n-2}{n}\Gamma_n)$. Moreover $N$ is a nef and big divisor.

Given a non-empty linear system $\ls$ which is not $(-1)$-special, we may always assume that it does 
not have negative product with any $(-1)$-curve and with $\pi^*\Gamma_n$. What we want to prove in this case is that 
$h^1(L)=0$. \\
Let $L$ be the line bundle associated to the proper transform of the system $\ls$ on
the blow-up surface $S$. If $C$ is an irreducible curve such that $L \cdot C<0$, then  $C$ is contained in the fixed part of $\ls$ and $C^2<0$. 
This means that $2g(C)-2=C^2+C\cdot K_S<0$, which implies that
$g(C)=0$ and $C^2$ equal $-1$ or $-2$. By hypothesis $C^2$ can not be a
$(-1)$-curve so that it has to be a $(-2)$-curve. From $0=C\cdot (N+\Delta)$ and the fact that $N$ is nef we deduce that $C\cdot\Delta=0$ so that $C\cdot (\pi^*(H+2F)-\sum_i E_i)=0$.
The last equality implies that $C\cdot (\pi^*(H+F)-\sum_i E_i)<0$ which is a contradiction since $\pi^*(H+F)-\sum_i E_i$ is a $(-1)$-curve, so that it is irreducible and distinct from $C$.

We proved that $L$ is nef so that $N+L$ is nef and big.  We conclude by observing that $L=K_S+(L-K_S)=K_S+(L+N)+\Delta$, so that $h^1(L)=0$, by KawamataÐViehweg vanishing theorem (see~\cite[Theorem 9.1.18]{laz}).
\end{proof}

\section{$(-1)$-Curves}
In what follow, with abuse of language, we call {\em $(-1)$-curve} also
the curves of $\f_n$ whose strict transform is a $(-1)$-curve on the blow
up surface $S$.
In order for $E\in \ls_n(a,b,m^r)$ to be a $(-1)$-curve, we must have:
\begin{equation}
\label{sistem1}
\left\{ \begin{array}{l}
        E^2=-1 \\
        E\cdot K=-1 \\  
        \end{array}
\right.  
\end{equation}  
which leads to the system: 
\begin{equation}
\label{sistem2}
\left\{ \begin{array}{l}
        2ab+nb^2-rm^2+1=0 \\
        2a+nb+2b-rm-1=0 \\  
        \end{array}
\right.  
\end{equation}  

By eliminating $a$, we obtain 
\begin{equation}
\label{equ1}
2b^2-rmb-b+rm^2-1=0.
\end{equation}
This is the equation of an irreducible conic in the coordinates $b,r$ for
$m>1$.
For $m>1$, observe that $r=0,\ b=1$ is an integral solution of
(\ref{equ1}), hence
we obtain a rational parameterization which gives us all the rational solutions 
of (\ref{equ1}). Let $b=rt+1$ where $t=p/q$ with $p,q\in\z$ and 
$\gcd(p,q)=1$, substituting in (\ref{equ1}) we obtain:
\[
        b=\frac{m^2p+q}{mp-2q} \hspace{1cm}
        r=\frac{p}{q}\ \frac{p(m^2-m)+3q}{pm-2q}.    
\]  
Now we look for integral solutions,
from the expression of $r$ it follows that $q\ |\ m(m-1)$ because $p$ and
$q$ are coprime and $q\ |\ p(m^2-m)+3q$. This allows us to calculate all the 
possible values of $q$ and also all the integral solutions of
(\ref{sistem2}) for a given $m$.\\ 
For $m=1$, $2b^2-rb-b+r-1=(b-1)(2b+1-r)$. $b=1$ gives the solution:
$\ls_n(e,1,1^{2e+n+1})$. 
If $r=2b+1$ substituting in the second equation of (\ref{sistem2}) we
obtain
$a=(2-nb)/2$, which is negative unless one has either $b\leq2$ or $n=0$.
If $n=0$ then $a=1$ and the system is $\ls_0(1,b,1^{2b+1})$ which is the
same as
$\ls_0(b,1,1^{2b+1})$ since on $\f_0$, $F$ and $H$ have the same properties. so
this system is already contained in $\ls_n(e,1,1^{2e+n+1})$.
If $n>0$ the systems are: $\ls_n(1,0,1)$, $\ls_2(0,1,1^3)$, $\ls_1(0,2,1^5)$.
Also $\ls_2(0,1,1^3)$ is contained in $\ls_n(e,1,1^{2e+n+1})$.
These are all the homogeneous $(-1)$-curves with multiplicity one: \\

\begin{center}
\fbox{
$\ls_1(0,2,1^5)\hspace{1cm} \ls_n(1,0,1)\hspace{1cm} \ls_n(e,1,1^{2e+n+1})$}
\end{center}

\section{Elementary Transformations}

In this section we consider a well known class of birational transformations
from an $\f_n$ surface to an $\f_{n-1}$: the {\em elementary transformation}.
These transformations have two good properties for our scope:

\begin{itemize}
\item The virtual dimension of linear systems are preserved through any such
transformation.
\item The $(-1)$-curves goes into $(-1)$-curves.
\end{itemize}

Sometimes one has also the following:
\begin{itemize}
\item The effective dimension of linear systems is preserved.
\end{itemize}

These transformations, allow us in some cases, to translate the problem of
homogeneous linear systems on $\f_n$ to a problem of quasi\ -\ homogeneous linear
systems on $\p^2$.\\

In order to understand the action of such transformations on linear systems, 
consider a surface $\f_n$ and a point $p$ not on the $(-n)$-curve. 
Blow up the fiber through $p$ and then blow down the strict transform 
of the fiber. 
Consider the system: $\ls_n(a,b,m)$ with only one singular point $p$. Let
$\pi_1$ be the blow up of $p$ with exceptional divisor $E$. 
Let $H_n,\ F_n,\ H_{n-1},\ F_{n-1}$ be the generators of $\pic(\f_n)$ and
$\pic(\f_{n-1})$ respectively, with $H_i^2=i,\ H_i\cdot F_i=1,\ F_i^2=0$
with $i\in\{n-1,n\}$. Let $F$ be the
strict transform of the fiber  $F_n$ through the point $p$ and let
$H=\pi_1^*H_n$. Now let $\pi_2$ be the blow down of $F$, we obtain:
\[
\pi_1^*H_n=H,\  \pi_1^*F_n=E+F,\  \pi_2^*H_{n-1}=H-E,\
\pi_2^*F_{n-1}=E+F.
\]  
The strict transform of a curve belonging to the system $\ls_n(a,b,m)$ is
$a\pi_1^*F_n+b\pi_1^*H_n-mE=a(F+E)+bH-mE=(a-m+b)(E+F)+b(H-E)-(b-m)F=
(a-m+b)\pi_2^*F_{n-1}+b\pi_2^*H_{n-1}-(b-m)F$. 
After blowing down it becomes $\ls_{n-1}(a-m+b,b,b-m)$.
Observe that if $b>m$ then in the last system the point of multiplicity $b-m$ 
belongs to the $1-n$ curve.
This gives the following transformation for $k<n$ points:
\begin{equation} 
\label{elemtransf}
\ls_n(a,b,m^r)\rt \ls_{n-k}(a+k(b-m),b,m^{n-k},(b-m)^k)
\end{equation}
It is a simple calculation to verify that the virtual dimension of such systems
is the same. If $C$ is a $(-1)$-curve of $\f_n$ through $r$ points, 
after an elementary transformation (with center in one of 
the $r$ points), the self-intersection and the genus of the new curve through $r-1$ points
are the same, hence these transformations preserve $(-1)$-curves. \\

Now we are able to prove the conjecture in the following cases:
\begin{propo}
\label{propo1}
The system $\ls_n(a,b,m^r)$ with $b\leq m+1$ is special if and only if it
is $(-1)$-special.
\end{propo}

\begin{proof}
By \ref{propo:<=n+1} we need to consider only systems with more than $n+1$ points.
We distinguish three cases: \\

$\bullet\ {b<m}$. In this case, since
$\ls_n(a,b,m^r)\cdot \ls_n(1,0,1)=b-m<0$, then the fibers $F$ through the
singular points are fixed components of $\ls_n(a,b,m^r)$. 
Let $k=m-b$, the system consists of a fixed part and of the residual
system $\ls_n(a-kr,m-k,(m-k)^r)$. 
If $b\leq m-2$ and the residual system is not empty, then the system is
$(-1)$-special. If $b=m-1$ then the residual system is $\ls_n(a-r,m-1,(m-1)^r)$.
In the last system we have $b=m-k$, so we need only to study the case $b=m$.\\
$\bullet\ {b=m}$. By making a sequence of $n-1$ elementary transformations we obtain the 
system $\ls_1(a,m,m^{r-n+1})$. Blowing down the $(-1)$-curve of $\f_1$ we obtain the 
quasi-homogeneous system on $\p^2:\ \ls_{\p^2}(a+m,\- a,\- m^{r-n+1})$, which
is special if and only if it is $(-1)$-special (see \cite[Proposition 6.2]{cm})\\
$\bullet\ {b=m+1}$.
Consider the system $\ls_n(a,m+1,m^r)$, after $n-1$ elementary transformations,
the system becomes $\ls_1(a+n-1,\- m+1,\- m^{r-n+1},\- 1^{n-1})$, with the $n-1$ points
belonging to the $(-1)$-curve. This system may be again transformed by blowing
up two points outside the $(-1)$-curve and blowing down the two fibers through
these points. In this way we obtain the system
$\ls_1(a+n-m,m+1,m^{r-n-1},1^{n+1})$ with the $n+1$ points belonging to a
section of self intersection $1$.
Blowing down the $(-1)$-curve we obtain the plane system
$\ls_{\p^2}(a+n+1,a+n-m,m^{r-n-1},1^{n+1})$
where the $n+1$ points belong to a line $R$.
Now we distinguish two cases, according to the fact that the system 
$\lss=\ls_{\p^2}(a+n+1,a+n-m,m^{r-n-1})$ is special or not.

If $\lss$ is special then, as shown in \cite[Proposition 6.4]{cm}, it is 
$(-1)$-special. This means that the initial system is $(-1)$-special. \\
If $\lss$ is non-special, but $\ls_{\p^2}(a+n+1,a+n-m,m^{r-n-1},1^{n+1})$ is special, 
this means that this last system contains
the line $R$ through the $q_i$'s points, otherwise these (generic) points of the
line, will impose independent conditions on the curves of $\lss$. 
This means than $h^0(\lss-\sum_{i=0}^nq_i)=h^0(\lss)$, hence we have the
following exact sequence:
\[
0\rt H^0(\lss_{\mid R})\rt H^1(\lss-R)\rt H^1(\lss)\rt H^1(\lss_{\mid R})\rt\cdots
\]
If $h^1(\lss-R)\neq 0$, then $\lss-R$ is special and hence $(-1)$-special 
(by \cite[Proposition 6.4]{cm}). This imply that also $\lss$ is $(-1)$-special.\\
If $h^1(\lss-R)=0$ then $h^0(\lss_{\mid R})=0$ and this happen iff
$\deg\lss_{\mid R}<0$,
which means that $n+1>\deg\lss$. This can not happen since $\deg\lss=a+n+1$.
\end{proof} 

\section{$(-1)$-Special Systems}
 
A $(-1)$-special homogeneous system $\ls$ may contain a $(-1)$-curve $A$
which is not homogeneous. In this case it must contain also all the 
$(-1)$-curves $A_i$ obtained from $A$ by a permutation of the base points. 
No two of the $A_i$ can meet, because if two $(-1)$-curves on a rational 
surface meet, their union moves in a linear system, and so the union cannot be 
part of the fixed divisor of $\ls$. 
Let $p_1,\cdots,p_r$ be the base points of the linear system,
since the Picard group of $S$ has rank $r+2$ there can be at most $r+1$ 
of these disjoint $(-1)$-curves. On the other hand let 
$k_i=\#\{p_i\mid \mu_{p_i}(A)=h_i\}$, where $\mu_{p_i}(A)$ is the multiplicity
of $A$ at the point $p_i$ which may be also $0$. In this way we have 
$\sum_{i=1}^s k_i=r$ and the number of distinct $A_i$ through the $r$ points 
is: $r!/(k_1!\cdots k_s!)$.
\begin{lemma}
Let $r,s,k_i\in\n$ such that $\sum_{i=1}^sk_i=r$, then $k_1!\cdots
k_s!\leq (r-s+1)!$ 
\end{lemma}

\begin{proof}
By hypothesis one has: $r-s+1=1+\sum_{i=1}^s(k_i-1)$. Define
$s_0=1$ and $s_t=1+\sum_{i=1}^t(k_i-1)$ for $t\geq 1$, observe that
\[
\prod_{j=1+s_t}^{s_{t+1}}j\geq k_{t+1}!
\]
because on the left side of the
inequality there are $k_{t+1}-1$ terms which are greater then or equal to
those on the right side (on this side the terms are $k_{t+1}$ but the
first is $1$ and does not give contribution in the product).
Taking products on $t$ we obtain the thesis.
\end{proof}
From the preceding lemma we obtain the following:
\begin{equation}
\label{compound}
\frac{r!}{k_1!\cdots k_s!}\geq\frac{r!}{(r-s+1)!}.
\end{equation}
The right side is a polynomial $P_s(r)$ in $r$ of degree $s-1$, it is easy
to prove that $P_s(r)>r+1$ for $r\geq s\geq 3$.  
If $s=1$ we have $P_1(r)=1$ and the system $A$ is homogeneous through the 
$r$ points. If $s=2$ then $P_2(r)=r$. 
In this case, the left hand of \ref{compound} is a binomial coefficient:
$\binom{r}{k_1}$. A simple calculation shows that $\binom{r}{2}\leq
r+1$ unless $r\leq 3$, but in this case $k_2\leq 1$. In the other cases we have 
$\binom{r}{k_1}\geq \binom{r}{2}$
for $2\leq k_1\leq [r/2]$. So we may assume that $k_1=1$ and $k_2=r-1$.
In this case the $(-1)$-curves are of type $\ls_n(a,b,m_1,m_2^{r-1})$.
Consider two of them with multiplicity $m_1$ at two different points, $p_i$ and
$p_k$ and call them $A_i$ and $A_k$ respectively.
The conditions $A_i\cdot A_k=0$ if $i\neq k$ and $A_i^2=-1$ give us 
$m_1=m_2\pm 1$. 
For our purposes we need only to consider those with $m_i\leq 1$ (since they
must be contained twice in the base locus of the linear system to give a $(-1)$
special system). So we have two possibilities: \\
$\bullet\ m_1=0$ and $m_2=1$, in this case the system must have multiplicity at least
$2r$, since it contains all these curves as fixed components twice. From $m\geq
2r$, if $m\leq 3$ we obtain $r=1$.\\
$\bullet\ m_1=1$ and $m_2=0$, in this case there are no restriction on the 
multiplicity of the system. We call such systems {\em compound}. \\

\begin{propo}
\label{classification}
All homogeneous $(-1)$-special systems with multiplicity $\leq 3$ on 
$\f_n$ are listed in the following table: 

\begin{table}[ht]
\caption{$(-1)$-special Systems with $m\leq 3$}
\begin{tabular}{ccc}\hline
System &  $v(\ls)$ & $l(\ls)$ \\ \hline
$\ls_1(0,4,2^5)$ & $-1$ & $0$ \\
$\ls_1(0,6,3^5)$ & $-3$ & $0$ \\
$\ls_5(1,4,3^{10})$ & $-1$ & $0$ \\
$\ls_6(0,4,3^{11})$ & $-1$ & $0$ \\
$\ls_n(2e,2,2^{2e+n+1})$ & $-1$ & $0$ \\
$\ls_n(e,0,2^r)$ & $e-3r$ & $e-2r$ \\
$\ls_n(4e+n+1,2,3^{2e+n+1})$ & $-1$ & $0$ \\
$\ls_n(3e+1,3,3^{2e+n+1})$ & $1$ & $2$ \\
$\ls_n(3e,3,3^{2e+n+1})$ & $-3$ & $0$ \\
$\ls_n(e,1,3^r)$ & $2e+n-6r+1$ & $2e+n-5r+1$ \\ 
$\ls_n(e,0,3^r)$ & $e-6r$ & $e-3r$ \\
\hline 
\end{tabular} 
\end{table}

\end{propo}

\begin{proof}
We consider homogeneous linear systems $\ls(a,b,m^r)$ which have negative products
with $(-1)$-curves. 
The first product is:
\[
\ls_1(a,b,m^r)\cdot\ls_1(0,2,1^5)<0,
\]
in this case the system is of the form $\ls_1(a,b,m^5)$, with $m=2,3$.
If $m=2$, then $2a+2b-10< 0$. The left hand of the preceding inequality is an even
number which can not be less of $-2$ (because the linear system $\ls_1(a,b,2^5)$ can not
contain the system $\ls_1(0,2,1^5)$ more than twice). This implies that 
$a+b=4$ and that $\ls_1(a,b,m^5)=2\ls_1(0,2,1^5)+\ls_1(a,b-4)$. Hence, in order for
the system to be not empty, it must be $b\geq 4$. This gives $b=4$ and $a=0$, so the
initial system is $\ls_1(0,4,2^5)$. \\ 
If $m=3$, then $\ls_1(a,b,3^5)\cdot \ls_1(0,2,1^5)\leq -1$ 
gives: $2(a+b)-15\leq -1$. 
If the product is equal to $-1$, then $a+b=7$ and the residual system is
$\ls_1(a,b-2,2^5)=\ls_1(a,5-a,2^5)$. This system is empty if $a\geq 3$ and
non-special if $a\geq 2$.
If the product is equal to $-3$, then $a+b=6$. The residual system 
is $\ls_1(a,b-6)$, hence $b\geq 6$, which implies that $b=6$ and $a=0$ and the 
system is $\ls_1(0,6,3^5)$. \\

The second product is:
\[
\ls_n(a,b,m^r)\cdot\ls_n(1,0,1)<0.
\]

If $m=2$ then, since the product is equal to $b-2$, $b$ may be equal to 
$0$ or $1$. 
\begin{itemize}
\item If $b=0$, then the residual system $\ls_n(a-2r,0)$ is not empty if and 
only if $a\geq 2r$. In this case 
the initial system is $(-1)$-special or empty. 
\item If $b=1$, then after removing the fixed part, the residual system 
$\ls_n(a-r,1,1^r)$ is non-special, since it has base points of multiplicity one. \\
In this case the initial system is non-special, since the $(-1)$-curves have
product equal to $-1$ with the initial system.
\end{itemize}

If $m=3$ then, since the product is equal to $b-3$, $b$ belong to the set 
$\{0,1,2\}$. 
\begin{itemize}
\item If $b=0$, then the residual
system $\ls_n(a-3r,0)$ is not empty if and only if $a\geq 3r$. In this case 
the initial system is $(-1)$-special or empty. 
\item If $b=1$, then the residual system $\ls_n(a-2r,1,1^r)$ is non-special of
virtual dimension $2a+n-5r+1$. Hence the initial system is $(-1)$-special or
empty.
\item If $b=2$, then the residual system is $\ls_n(a-r,2,2^r)$. 
This system may be $(-1)$-special, only if it has negative product with some of
the other $(-1)$-curves of the $\f_n$, hence it is already studied in the other 
cases of the present classification. 
\end{itemize}

The third product is:
\[
\ls_n(a,b,m^{2e+n+1})\cdot \ls_n(e,1,1^{2e+n+1})<0
\]
Here we consider only systems with $b\geq m$, since the other systems where
already considered in the preceding case.
We begin by considering the case $m=2$. \\

In this case, the product may be equal to $-1$ or to $-2$. 
In the first case one has $a=e(4-b)+n(2-b)+1$ and the 
residual system is $\ls_n(a-e,b-1,1^{2e+n+1})$.  
\begin{itemize}
\item If $b=2$ then $a=2e+1$ and the residual system is
$\ls_n(e+1,1,1^{2e+n+1})$. 
This system is non-special by proposition \ref{mult.one}.
\item If $b=3$ then $a=e-n+1$ and the residual system is
$\ls_n(-n+1,2,1^{2e+n+1})$. If $n\leq 1$, then by proposition \ref{mult.one}
the sistem is non-special. If $n\geq 2$, then the residual system has negative 
product with the $-n$-curve $\Gamma_n$ of $\f_n$, hence it is equal to
$\Gamma_n+\ls_n(1,\- 1,\- 1^{2e+n+1})$. By proposition \ref{mult.one}, this system 
is non special of dimension $2-2e$. This means that $e\leq 1$. Since
$a=e-n+1\geq 0$, the only possibility with $n\geq 2$ is $e=1$ and $n=2$. In this
case the initial system is $\ls_2(0,3,2^5)$, it has virtual dimension $0$, hence
it is non-special.
\item If $b=4$ then $a=-2n+1$ is negative unless $n=0$ and in this case the
residual system is empty.
\item If $b\geq 5$ then $a$ is always negative and the system is empty.
\end{itemize}
If the product is equal to $-2$, then $a=e(4-b)+n(2-b)$ and the residual system
is $\ls_n(a-2e,b-2,0)$. 
\begin{itemize}
\item If $b=2$ then $a=2e$ and the residual system $\ls_n(0,0,0)$ is not empty
and the system $\ls_n(2e,2,2^{2e+n+1})$ is $(-1)$-special.
\item If $b=3$ then $a=e-n$ and the residual system $\ls_n(-e-n,1,0)$ is empty.
\item If $b\geq 4$ then the system is always empty.
\end{itemize}

Now we consider the case $m=3$.\\

In this case the product may be equal to $-1, -2, -3$.
In the first case $a=e(6-b)+n(3-b)+2$ and the residual system is
$\ls_n(a-e,b-1,2^{2e+n+1})$. 
\begin{itemize}
\item If $b=3$ then $a=3e+2$ and the residual system is
$\ls_n(2e+2,2,2^{2e+n+1})$. By \ref{propo1} this system is non-special.
\item If $b=4$ then $a=2e-n+2$ and the residual system is
$\ls_n(e-n+2,3,2^{2e+n+1})$. If $e-n+2\geq 0$ then by \ref{propo1} the system is
non-special. If $e-n+2 <0$ then the system has negative product with the $-n$
curve $\Gamma_n$. So $\ls_n(e-n+2,3,2^{2e+n+1})=\Gamma_n +
\ls_n(e+2,2,2^{2e+n+1})$. Since 
$\ls_n(e+2,2,2^{2e+n+1})\cdot\ls_n(e,1,1^{2e+n+1})=-e$, then $e\leq 2$ otherwise
the system is empty.\\
If $e=0$ then $a=-n+2<0$ because we suppose that $e-n+2<0$,
hence there is no such a possibility. \\
If $e=1$ then $a=4-n$ and $n\geq 4$. This
means that $n=4$ and $a=0$. This gives the system $\ls_4(0,4,3^7)$.
This system is equal to
$\ls_4(1,1,1^7)+\ls_4(-1,3,2^7)=\ls_4(1,1,1^7)+\ls_4(-4,1,0)+\ls_4(3,2,2^7)$
Since $\ls_4(3,2,2^7)\cdot\ls_4(1,1,1^7)=-1$, the initial system is equal to
$2\ls_4(1,1,1^7)+\ls_4(-4,1,0)+\ls_4(2,1,1^7)$. By proposition \ref{mult.one},
the last residual system is non-special of dimension $2$. Since the virtual
dimension of the initial system is $2$, the system is non-special.\\
If $e=2$ then $a=6-n$ and $n\geq 5$. This means that $n$ may be equal to $5$ or
$6$. If $n=5$ then $\ls_5(1,4,3^{10})=\ls_5(2,1,1^{10})+\ls_5(-1,3,2^{10})$,
this is equal to $\ls_5(2,1,1^{10})+\ \ls_5(-5,1,0)+\ \ls_5(4,2,2^{10})$.  
Since $\ls_5(4,2,2^{10})\cdot\ls_5(2,1,1^{10})=-2$, we have that the initial
system is equal to $3\ls_5(2,1,1^{10})+\ls_5(-5,1,0)$. The virtual dimension of
the initial system is $-1$, hence the system is $(-1)$-special.
If $n=6$ then $\ls_6(0,4,3^{11})=\ls_6(2,1,1^{11})+\ls_6(-2,3,2^{11})$,
this is equal to $\ls_6(2,1,1^{11})+\ \ls_6(-6,1,0)+\ \ls_5(4,2,2^{11})$.  
Since $\ls_6(4,2,2^{11})\cdot\ls_6(2,1,1^{11})=-2$, we have that the initial
system is equal to $3\ls_6(2,1,1^{11})+\ls_6(-6,1,0)$. The virtual dimension of
the initial system is $-1$, hence the system is $(-1)$-special.
\item If $b=5$ then $a=e-2n+2$ and the residual system is
$\ls_n(-2n+2,4,2^{2e+n+1})$. This system has negative product with the
$-n$-curve only if $n\geq 2$, otherwise the system is non-special, since it does
not have negative product with any other $(-1)$-curve.
If $n=2$ then the residual system $\ls_2(-2,4,2^{2e+3})$ is equal to 
$\ls_2(-2,1,0)+\ls_2(0,3,2^{2e+3})$. The last residual system is not
$(-1)$-special of virtual dimension $6-6e$ and since $a=e-2\geq 0$, this
dimension is negative. Hence the initial system is not $(-1)$-special.
\item If $b=6$ then $a=-3n+2$ and the system is empty.
\end{itemize}
If the product is equal to $-2$, then $a=e(6-b)+n(3-b)+1$ and the residual 
system is $\ls_n(a-2e,b-2,1^{2e+n+1})$. 
\begin{itemize}
\item If $b=3$ then $a=3e+1$
\item If $b=4$ then $a=2e-n+1$ and the residual system is
$\ls_n(-n+1,2,1^{2e+n+1})$. This is equal to $\Gamma_n + \ls_n(1,1,1^{2e+n+1})$.
The last system has dimension $2-2e$, hence it is not empty only if $e=0,1$.  
\item If $b=5$ then $a=e-2n+1$ and the residual system is
$\ls_n(-e-2n+1,2,1^{2e+n+1})$.
\end{itemize}
If the product is equal to $-3$, then $a=e(6-b)+n(3-b)$ and the residual 
system is $\ls_n(a-3e,b-3,0)$. 
\begin{itemize}
\item If $b=3$ then $a=3e$ and the system
$\ls_n(3e,3,3^{2e+n+1})=3\ls_n(e,1,\- 1^{2e+n+1})$ is $(-1)$-special.
\item If $b=4$ then $a=2e-n$ and the residual system $\ls_n(-e-n,1,0)$ is empty.
\item If $b=5$ then $a=e-2n$ and the residual system $\ls_n(-2e-2n,2,0)$ is
empty.
\item If $b=6$ then $a=-3n$ is negative and there are no such systems. 
\end{itemize}

\end{proof}

\section{Degeneration of $\f_n$ Surfaces}

Let $\Delta$ be: $\{ z\in\co | \mid z\mid <1\}$, consider the product with
$\f_n$ and the two projections maps $\pi_i$. 
Now blowing up the rational normal curve $H$ contained in $\pi_1^{-1}(0)$ 
we obtain a threefold $X$ and two maps: 

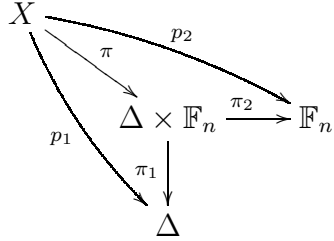
\[
\xymatrix{
X \ar@/_/[ddr]_{p_1} \ar@/^/[drr]^{p_2} \ar[dr]^{\pi} \\
& {}\Delta\times{\f_n} \ar[d]_{\pi_1}  \ar[r]^{\pi_2} & {}\f_n \\
& {}\Delta  &  }
\]

Let $X_t=p_1^{-1}(t)$, which for $t\neq 0$ is $\f_n$, while $X_0$ is 
the union of the proper transforms $\tilde{\f}$ of $\f_n$ and of the exceptional 
divisor $\f$ which is $\p(N_{\f_n \times\Delta\mid H})$. 
But $N_{\f_n \times\Delta\mid H}=\oc_{\p^1}(n) \oplus \oc_{\p^1}$. Hence the 
exceptional divisor is also an $\f_n$ that intersects the proper transform in a 
curve $\re$ whose class is $H$. Let us denote by $\rn$ the curve $R$ when
we consider it as a divisor on $\f$. To understand the class of $\rn$ in
$\pic (\f)$ observe that $(\f+\tilde{\f})\cdot\f=X_0\cdot\f=X_t\cdot\f=0$.
Hence $\f\cdot\f=-\tilde{\f}\cdot\f=-\rn$. 
Now consider the line bundle $\f$ and the two restrictions: 
${\f}_{\mid\f}=-\rn$ and ${\f}_{\mid\tilde{\f}}=\re$. If we consider the
restrictions of these systems to $R$ they must agree, since they come from 
the same line bundle. So is $-\rn^2=\re^2=n$, which tells us that $\rn$
is the $(-n)$-curve of $\f$. \\

Now consider a linear system $\mid aF+bH\mid$ on the generic fiber $X_t$, 
we have: 
$p_2^*(aF+bH)_{|\tilde{\f}}=aF+bH$ and $p_2^*(aF+bH)_{|\f}=(a+nb)F$.
So this is not sufficient to generate $\pic(\f)$.
If we consider the linear system
\[
{\mathcal{X}} (a,b,k)=\mid p_2^*(a F+b H)+k\tilde{\f}\mid,
\]
this restricted to $X_t$ gives again $\mid aF+bH\mid$, while restricted to 
$X_0$ gives $\mid aF+(b-k)H\mid$ on $\tilde{\f}$ and $\mid (a+nb-nk)F+kH\mid$ 
on $\f$. 
Now consider an homogeneous linear system $\ls=\ls_n(a,b,m^r)$.
We define a {\em $(k,s)$-degeneration} by sending $s$ of the $r$ points
to $\tilde{\f}$ and considering the system ${\mathcal{X}}(a,b,k)_{\mid X_0}$ with
$s$ base points on $\tilde{\f}$ and $r-s$ base points on $\f$.
Now let ${\mathcal{X}} = p_2^*\ls$, since the map $p_1$ is a flat morphism, by
semi-continuity is: 
\[
h^0(\ls)=h^0(X_t,\pi_2^*\ls)=h^0(X_t,{\mathcal{X}})\leq
h^0(X_0,{\mathcal{X}}).
\]
We call $\ls_0={\mathcal{X}}_{\mid X^0}$ and let $l_0=\dim\ls_0$.
This is the key of the method: in order to prove the conjecture for a
supposed non special system, one proves that $v(\ls)=l_0$ and this implies 
that also $v(\ls)=l(\ls)$. \\

In order to give a linear system on $X_0$, one has to give two linear system on
$\f$ and $\tilde{\f}$ which agree on $R=\f\cap\tilde{\f}$. 
If we call $\rfn$ (respectively $\rfe$) the restrictions of $\lfe$ (respectively
$\lfn$) to $R$, then we have the
following exact sequences:
\[
0\rt\klfe\rt\lfe\rt\rfn\rt 0
\]
\[
0\rt\klfn\rt\lfn\rt\rfe\rt 0.
\]
Where $\klfe$ and $\klfn$ are the kernels of the restrictions to $R$.
So we obtain the four systems:
\begin{center}
\fbox{
\parbox[c]{1cm}{
\begin{tabbing}
$\lfe = \ls_n(a,b-t,m^{r-s})$\hspace{1cm} \= $\lfn = \ls_n(a+n(b-t),t, m^s)$ \\ 
$\klfe = \ls_n(a,b-t-1,m^{r-s})$ \> $\klfn = \ls_n(a+n(b-t+1),t-1,m^s)$ 
\end{tabbing}
}}
\end{center}
\vspace{.3cm}

Now, in order to evaluate $h^0(X_0,\ls)$ we must consider the fibered
product:
\[
H^0(X_0,\ls)=H^0(\f,\lfn)\otimes_{H^0(R,\ls)} H^0(\tilde{\f},\lfe),
\]
which is obtained by taking the sections of $H^0(\f,\lfn)$ and 
$H^0(\tilde{\f},\lfe)$ whose restrictions to $R$ give the same element of 
$H^0(R,\ls)$.
Hence one has:
\[
l_0=\dim (\rfn\cap\rfe)+\dklfn +\dklfe +2. 
\]
In order to evaluate $\dim (\rfe\cap\rfn)$,
we recall the following lemma proved in \cite[proposition (6.1)]{cm}.
\begin{lemma}
Let $L\in\pic(\p^1)$, given two vector subspaces:
$V_1,V_2\subset H^0(\p^1,L)$ which are not transverse, there exists an
isomorphism $\psi$ of $\p^1$ such that $\psi^*V_1$ is transverse with
$V_2$.
\end{lemma}

Now, observe that such an isomorphism $\psi$ of $\p^1$ may be extended to
an isomorphism $\tau: \f_n\rt\f_n$ such that $\pi\circ\tau=\psi\circ\pi$
where $\pi: \f_n\rt\p^1$ is the projection map. To see this, it is sufficient
to consider the isomorphism $\psi$ as defined between two rational section
of the $\f_n$ and to extend it with the identity map on the fibers.

\begin{remark}
\label{cond1}
Let $v$ be the virtual dimension of the system $\ls_n(a,b,m^r)$.
Let $\drfe = \dlfe - \dklfe -1$ and $\drfn = \dlfn - \dklfn -1$
be the dimensions of the restrictions (to $R$) of the linear systems 
$\lfe$ and $\lfn$ respectively. Let $\vfe,\ \vfn,\ \kvfe,\ \kvfn$ be the virtual
dimensions of the systems $\lfe,\ \lfn,\ \klfe,\ \klfn$. 
An easy calculation shows that the following identities hold:
\begin{itemize}
\item[(i)] $\vfe + \vfn = v+a+n(b-k)$
\item[(ii)] $\kvfn + \vfe=\kvfe + \vfn = v - 1$ 
\item[(iii)] If $\drfe+\drfn\leq a+n(b-k)-1$, then $l_0=\dklfe +\dklfn
+1$.
\item[(iv)] If $\drfe+\drfn\geq a+n(b-k)-1$, then $l_0 =\dlfe+\dlfn
-a-n(b-k)$.
\end{itemize}
\end{remark}

\section{The Main Theorem}

Now we are ready to prove the main conjecture in the case $m=2,3$.
We need first a lemma:
 
\begin{lemma}
\label{lemma2}
For each $r\leq n+1$, the linear system $\ls_n(a,b,m^r)$ is special 
if and only if it is $(-1)$-special.
\end{lemma}

\begin{proof}
This is due to the fact that $K_S$ is nef. 
\end{proof}  

\begin{teo}
\label{mainteo}
Every special homogeneous system of multiplicity $\leq 3$ on an $\f_n$ 
surface, is a $(-1)$-special system.
\end{teo}   

\begin{proof} 
\begin{center}
Multiplicity $2$
\end{center}
The case $n=1$ is the same as the case of quasi homogeneous (i.e. systems
with all except one points of the same multiplicity) systems on $\p^2$,
which is completely classified in \cite{cm} for multiplicity $\leq 3$. 
For the remaining $n$, we proceed by induction on the number $r$ of
points.
For systems with only one point we proved the thesis in \ref{lemma2}. 
By lemma \ref{propo1} we may consider the linear system $\ls=\ls_n(a,b,2^r)$
with $b\geq 4$. 
First consider the case $v(\ls_n(a,b,2^r))<0$.
Performing a $(1,s)$-degeneration, we obtain the four systems: 
\[
\lfe = \ls_n(a,b-1,2^{r-s}),\ \lfn = \ls_n(a+n(b-1),1, 2^s)
\]
\[
\klfe = \ls_n(a,b-2,2^{r-s}),\ \klfn = \ls_n(a+nb,0,2^s).
\] 
{\em Step I}. We want to find an $s$ such that
\begin{equation}
\label{step1}
\klfn=\klfe=\emptyset.
\end{equation}

Observe that $\klfn$ is a $(-1)$-special system of dimension $a+nb-2s$, so
it is empty iff $s>(a+nb)/2$. There always exists such an $s$ 
because, from the inequalities $v(\ls)<0$ and $b\geq 3$ it follows that
$r>(a+nb)/2$.
A sufficient condition to have $\klfe=\emptyset$ is that it is
non special and $v(\klfe)\leq v(\ls)$. This gives 
$s\leq 2(a+nb-n+1)/3$. So $s$ must satisfy:

\begin{equation}
\label{equa}
\frac{a+nb}{2}<s\leq \frac{2}{3}(a+nb-n+1)
\end{equation}

A sufficient condition for $s$ to exist is that
$2(a+nb-n+1)/3-(a+nb)/2\geq 1$. Solving the inequality we
obtain: $n(b-4)\geq 2-a$ which is always true if $b\geq 5$ or if $b=4$ and
$a\geq 2$. In the remaining case, $b=4,\ a=1$ it is again possible to find 
$s=2n+1$. 
$\klfe$ may be  $(-1)$-special in only one case (see Proposition
\ref{classification}) with $b\geq 4$: $\klfe=\ls_n(2e,2,2^{2e+n+1})$. 
For $e\geq 4$ there are two values of $s$ that satisfy (\ref{equa}), 
so we may choose the value for which $\klfe$ is
non special. In the remaining cases ($e\leq 3$) $v(\ls)\geq 0$. 
Now we have chosen $s$ such that $\klfe$ and $\klfn$ are empty and we
finish step I. \\

{\em Step II}. We want to prove that $\drfe+\drfn\leq a+n(b-1)-1$
which by [\ref{cond1},iii] ends the proof.
Observe that (Proposition \ref{classification}) $\lfn$ is
non special $\lfe$ is non special for $b\geq 4$. 
Then by [\ref{cond1},i], $\dlfe+\dlfn=a+n(b-1)+v$ and since $v\leq -1$, by
[\ref{cond1},iii] we finish.\\

Suppose now that $v(\ls_n(a,b,2^r))\geq 0$ and consider
again a $(1,s)$-\ degeneration. In this case we need the regularity
of $\lfn$ and $\lfe$ and the inequality:
\[
\drfe +\drfn\geq a+n(b-1)-1.
\]
So applying [\ref{cond1},i] and [\ref{cond1},iv] we prove that $l_0=v$ and
we finish.
Choosing $s\leq (a+nb)/2$, then both $\lfn$ and $\klfn$ are not empty
and by proposition \ref{classification} for $b\geq 4$, $\lfn$ and $\lfe$
are non special. 
$\klfn$ is $(-1)$-special of dimension $a+nb-2s$. 
and $\dlfn=2(a+bn)-n-3s+1$ and 
This gives $\drfn = \dlfn - \dklfn - 1 = a+bn-n-s$. So we have 
only to check that $\drfe\geq s-1$. If $\klfe$ is non special the inequality is 
true for $b\geq 2$. 
If $\klfe$ is $(-1)$-special the only possibility with $b\geq 4$ is:  
$\ls_n(a,b-2,2^{r-s})=\ls_n(2e,2,2^{2e+n+1})$ which has dimension $0$.
In this case $\lfe$ is non special of dimension $2e+3n$. So $\drfe\geq n+s-1$
iff $s\leq 2e+2n$ which is always true.

\begin{center}
Multiplicity $3$
\end{center}

The technique is the same as in the multiplicity $2$ case, so we proceed
again by induction assuming that $b\geq 5$ by lemma \ref{propo1}.
First suppose that $v(\ls_n(a,b,3^r))<0$. 
In this case, we make a $(2,s)$-degeneration, obtaining the four systems: 
\[
\lfe = \ls_n(a,b-2,3^{r-s}),\ \lfn = \ls_n(a+n(b-2),2,3^s)
\]
\[
\klfe = \ls_n(a,b-3,3^{r-s}),\ \klfn = \ls_n(a+n(b-1),1,3^s)
\]

{\em Step I}. We need to find an $s$ that satisfies (\ref{step1}).
Observe that by \ref{classification} $\klfn$ is a $(-1)$-special system of
dimension $2(a+nb)-n-5s+1$. 
So choosing $s>(2(a+nb)-n+1)/5$, we obtain: $\klfn=\emptyset$.
It is possible to chose such an $s$ since $v<0$ and $b\geq 5$ imply that
$r>(2a+nb+2)/2-1/6$.
Now, solving the inequality $v(\klfe)\leq v(\ls)$ 
we obtain $s\leq (a+bn-n+1)/2$. So we want:
\begin{equation}
\label{equ5}
\frac{2(a+nb)-n+1}{5}<s\leq\frac{a+bn-n+1}{2}, 
\end{equation}
A sufficient condition for $s$ to exist is that
$(a+bn-n+1)/2-(2(a+nb)-n+1)/5\geq 1$, which is equivalent to 
$(b-3)n+a-7\geq 0$. This inequality is not satisfied only if $b=5,\ n=2$ and
$a\leq 2$, but in these cases $s=5$ works.
If $\klfe$ is $(-1)$-special, then by \ref{classification} it is one of
the following: $\ls_n(4e+n+1,2,3^{2e+n+1}),\ 
\ls_n(3e+1,3,3^{2e+n+1}),\  \ls_n(3e,3,3^{2e+n+1})$. \\
 
First we try to solve in each case the inequality:
$(a+bn-n+1)/2-(2(a+nb)-n+1)/5\geq 2$. This allow us to chose $s$
in
two ways and to avoid the $(-1)$-special case. 
In the first case we obtain $4e+3n\geq 16$. This is true if $e\geq 3$, in the
remaining cases, making the explicit calculation, there are two values for $s$. 
In the second case we have $3e+3n\geq 16$, which is true if $e\geq 4$. If 
$e\leq 3$ it is again possible to find two values for $s$ except for $e=1,\ n=2,
\ s=7$ but in this case $v(\ls)=4$.
In the third case we have $3e+3n\geq 17$ which is always true for $e\geq 4$. In
the remaining cases the only exceptions are $e=1,\ n=3,\ s=9$ and $e=2,\ n=2,\
s=8$. But in these cases $v(\ls)=0$.
So we may assume that also $\klfe=\emptyset$. \\

{\em Step II}. We want to prove that $\drfe+\drfn\leq a+n(b-1)-1$
which by [\ref{cond1},iii] ends the proof.
$\lfe$ is $(-1)$-special only if $\lfe=\ls_n(3e,3,3^{2e+n+1})$ or $\lfe=\ls_n(3e+1,3,3^{2e+n+1})$.
In the first case $\dlfe=0$ and $\lfn=\ls_n(3(e+n),2,3^s)$. If $\lfn$ is not
$(-1)$-special, then
$\vfn=9e+12n-6s+2$, hence $\dlfe+\dlfn=\vfn$. From (\ref{equ5}) we deduce that 
$\vfn<(9e+6n-2)/5<3e+3n=a+n(b-2)$, so the system is empty.
If $\lfn=\ls_n(3(e+n),2,3^s)$ is $(-1)$-special then it must have dimension $0$,
so $\drfe+\drfn = \dlfn+\dlfe=0<a+n(b-2)$ and the system is empty. \\
If $\lfe=\ls_n(3e+1,3,3^{2e+n+1})$, then $\dlfe=2$.
If $\lfn=\ls_n(3(e+n)+1,2,3^s)$ is non special, then $\vfn=9e+12n-6s+5$. 
From $s>(2(a+nb)-n+1)/5$ we deduce that 
$\dlfe+\dlfn=\vfn+2<(9e+6n+1)/5+2<3e+3n=a+n(b-2)$
and the system is empty.
If $\lfn=\ls_n(3(e+n)+1,2,3^s)$ is $(-1)$-special then it may have dimension $0$ or $2$ 
and again $\dlfn+\dlfe<a+n(b-2)$ and the system is empty. \\

So we may assume that $\lfe$ is non special. If also $\lfn$ is non special, then 
$\dlfn+\dlfe=\vfn+\vfe=a+n(b-2)+v(\ls)<a+n(b-2)$ and again the result
follows. \\

$\lfn$ may be $(-1)$-special only if $\lfn=\ls_n(4e+n+1,2,3^{2e+n+1})$. In this case
we have $a+n(b-2)=4e+n+1$ and $s=2e+n+1$. It is easy to prove that we may choose 
$s$ in two ways in order to avoid the $(-1)$-special case. \\

Now we study the case $v(\ls_n(a,b,3^r))\geq 0$.
We perform a $(3,s)$-degeneration, obtaining the systems:
\[
\lfe = \ls_n(a,b-3,3^{r-s}),\ \lfn = \ls_n(a+n(b-3),3,3^s)
\]
\[
\klfe = \ls_n(a,b-4,3^{r-s}),\ \klfn = \ls_n(a+n(b-2),2,3^s)
\]
Now let $\nu\equiv a+n(b-1)+1\pmod{2}$ and let $s=(a+n(b-1)+1+\nu)/2$.
With this choice of $s$ the system $\klfn$ is empty or $(-1)$-special with 
$\kvfn=-1$ and $\dklfn=0$.
The system $\lfn$ is always non special and $\vfn=a+n(b-3)-3\nu$. Then
$\drfn=a+n(b-3)-3\nu+\epsilon$ where $\epsilon=1$ if the system is $(-1)$-special,
otherwise is $0$. \\

Observe that $\vfe-v=3\nu\geq 0$, hence $\lfe$ is not empty.
If $\lfe$ is non special and $\klfe=\emptyset$ then $\drfe=v+3\nu$ and
$\drfe+\drfn=a+n(b-3)+\epsilon$ so $l_0 = \dlfe+ \dlfn -a-n(b-3)=v$.
If $\klfe\neq\emptyset$ is non special, then $\drfe=\vfe-\kvfe-1=a+n(b-3)$ and
again we finish.
If $\klfe\neq\emptyset$ is $(-1)$-special then for all the possibilities we have 
$\dklfe-\kvfe\leq 5$. Hence $\drfe\geq a+n(b-3)-5$ which gives 
$\drfe+\drfn\geq a+n(b-3)-1$ which is true if $a+n(b-3)\geq 8$. \\

If $\lfe$ is $(-1)$-special then it may be only $\lfe=\ls_n(3e+1,3,3^{2e+n+1})$.
In this case choosing $s=(a+n(b-1)+1+\nu)/2+1$ we have $\klfn=\emptyset$
and $\lfn=\ls_n(3e+4n+1,2,3^{s+1})$ which is non special. Hence 
$\drfn=\vfn=3(e+n-\nu)-5$. $\klfe=\ls_n(3e+n+1,2,3^{2e+n})$ is non special and empty. 
$\lfe=\ls_n(3e+1,3,3^{2e+n})$ is non special and $\vfe=7$. Therefore 
$\drfe+\drfn=3(e+n-\nu)+2$, hence if $\nu=0$ we finish. If $\nu=1$ then
$\ls=\ls_n(3e+1,6,3^r)$ with $r=(7e+7n+5)/2$ and in this case $v=-2$, a
contradiction.
\end{proof}

\bibliographystyle{plain}

\end{document}